\documentclass[12pt]{article}
\usepackage{tikz}
\usetikzlibrary{intersections,calc,arrows}
\usepackage{amssymb}
\date{}
\makeatletter
\newcommand{\figcaption}[1]{\def\@captype{figure}\caption{#1}}
\newcommand{\tblcaption}[1]{\def\@captype{table}\caption{#1}}

\@addtoreset{equation}{section}
\makeatother
\newcommand{\qed}{\hbox{\rule[-2pt]{3pt}{6pt}}}

\begin{document}
\title {\bf Exact solutions and bifurcation curves of nonlocal elliptic equations 
with convolutional Kirchhoff functions}

\author{{\bf Tetsutaro Shibata}
\\
{\small Laboratory of Mathematics, 
Graduate School of Advanced Science and Engineering,} 
\\
{\small Hiroshima University, Higashi-Hiroshima, 739-8527, Japan}
}

\maketitle
\footnote[0]{E-mail: tshibata@hiroshima-u.ac.jp}

\vspace{-0.5cm}

\begin{abstract}
We study the one-dimensional nonlocal elliptic equation of Kirchhoff type with 
convolutional Kirchhoff functions. 
We establish the exact  solutions $u_\lambda$ and  bifurcation curves $\lambda(\alpha)$, 
where $\alpha:= \Vert u_\lambda\Vert_\infty$.

\end{abstract}


\vspace{0.5cm}


\section{Introduction} 		      

We consider the following one-dimensional nonlocal elliptic equation with 
convolutional Kirchhoff function

\begin{equation}
\left\{
\begin{array}{l}
-\left(\displaystyle{\int_0^1}  f(x)u(x)^q dx\right) u''(x)= \lambda u(x)^p, \enskip 
x \in I:= (0,1),
\vspace{0.1cm}
\\
u(x) > 0, \enskip x\in I, 
\vspace{0.1cm}
\\
u(0) = u(1) = 0,
\end{array}
\right.
\end{equation}
where $f(x) = (1-x)^n$ ($n \in \mathbb{N}$) and $p, q > 1$ are given constants. 
Further, 
$\lambda > 0$ is a bifurcation parameter.

\noindent
Equation (1.1) is motivated by the convolution nonlocal elliptic problem of Kirchhoff type in [9]:
\begin{equation}
-J\left((h *u^q)(1)\right)u''(x)= 
\lambda g(x, u(x)), \enskip 
x \in I,
\end{equation}
where $J$, $h$ and $g$ are continuous functions and 
$$
(h * u^q)(t) := \int_0^t h(t-s)u(s)^qds.
$$
If we put $J(y) = y$, $h(y) = y^n$ and $g(x, u) = u^p$, then we obtain (1.1). 

The purpose of this paper is to obtain the exact solutions and bifurcation curves of the problem (1.1) by concentrating 
on the typical convolutional equation (1.1). Our results are novel, since it seems few results 
to treat such problem as (1.1) from a view point of bifurcation analysis, and the results 
obtained here will be the good first step to understand well the structures of the solutions and  
bifurcation curves in the field of non-local elliptic problems. Moreover, as far as the author knows, the bifurcation phenomena of 
non-local problem with the coefficient 
coming from the convolution has not been considered before.

It is known that there are so many results concerning non-local problems. We refer to 
[1--4, 6--11, 13--15, 20] and 
the references therein. Add to this, 
there are many interesting and significant 
motivation to study this kind of non-local problems as (1.2). We refer to [9] and the references therein to 
understand the background of this problem. On the other hand, although the 
analysis of bifurcation diagrams  
are very popular problem in nonlinear elliptic problems, it seems a few results 
concerning the bifurcation problems for non-local problems. We refer to [16-19, 21].

Before stating our results, we explain some notations. 
For $p> 1$, let 
\begin{equation}
\left\{
\begin{array}{l}
-W''(x)= W(x)^p, \enskip 
x \in I,
\vspace{0.1cm}
\\
W(x) > 0, \enskip x\in I, 
\vspace{0.1cm}
\\
W(0) = W(1) = 0.
\end{array}
\right.
\end{equation}
We know from [5] that there exists a unique solution $W_p(x)$ of (1.3). 
For $d, k\ge 0$, we put 
\begin{eqnarray}
L_{k,d}&:=& \int_0^1 \frac{s^d}{\sqrt{1-s^{k+1}}}ds,
\\
M_{k,d}&:=& \int_0^1 (1-x)^kW_p(x)^ddx, 
\\
R_{k, d}&:=& \int_0^1 x^k W_p(x)^ddx,
\\
S_{k, d}&:=& \int_0^{1/2} x^k W_p(x)^ddx.
\end{eqnarray}
Let $\Vert \cdot\Vert_\beta$ ($1 \le \beta \le \infty$) be the usual $L^\beta$-norm. 
We know from [19] that for $p > 1$, 
\begin{eqnarray}
\xi_p := \Vert W_p\Vert_\infty &=&  (2(p+1))^{1/(p-1)}L_{p,0}^{2/(p-1)}.
\end{eqnarray} 
Now we state our results.

\vspace{0.2cm}

\noindent
{\bf Theorem 1.1.} {\it Let $f(x) = (1-x)^n$ $(n \in \mathbb{N})$. Further, 
$\lambda > 0$ is a given constant. 

\noindent
(i) Assume that $p \not= q+1$. Then the solution $u_\lambda$ is given by 
\begin{eqnarray}
u_\lambda = \left(\frac{\lambda}{M_{n,q}}\right)^{1/(q-p+1)}W_p(x).
\end{eqnarray}

\noindent
{(ii)} Assume that $p = q+1$. 

\noindent
(a) Suppose that $\lambda = M_{n,q}$. Then all the solutions $u_\lambda$ 
are represented as $u_\lambda = tW_p$, where $t > 0$ is an arbitrary constant.

\noindent
{(b)} Assume that there exists a solution $U_\lambda$ of (1.1). Then 
$U_\lambda = r_\lambda W_p$, where $r_\lambda:= (Q_{n,q}/\lambda)^{1/(p-1)}$ and 
\begin{eqnarray}
Q_{n,q}:= \int_0^1 (1-x)^nu_\lambda(x)^q dx.
\end{eqnarray}
Moreover, $\lambda = M_{n,q}$ holds. Therefore, by (i) above, all the solutions $u_\lambda$ 
of (1.1) are 
obtained as $u_\lambda = tW_p$, where $t > 0$ is an arbitrary constant. 

\noindent
{(c)} Assume that $\lambda \not= M_{n,q}$. Then (1.1) has no solutions.

}

\vspace{0.2cm}

By Theorem 1.1, we see that the essential point to obtain the solution $u_\lambda$ is 
to find $M_{n,q}$. In the following Theorem 1.2 and 1.3, by using Theorem 1.1 (i), we
 obtain the exact solution $u_\lambda$ for given $\lambda > 0$ and show that $\lambda$ is parameterized by 
$\alpha:= \Vert u_\lambda\Vert_\infty$, namely, $\lambda =\lambda(\alpha)$ and 
establish the exact formula of $\lambda(\alpha)$.  

\vspace{0.2cm}

\noindent
{\bf Theorem 1.2.} {\it Let $f(x) = 1-x$ in (1.1). Assume that $p \not= q+1$. 
Then for any given $\lambda > 0$, 
\begin{eqnarray}
u_\lambda(x) 
&=& \lambda^{1/(q-p+1)}
\\
&&\times\left\{
2^{(q-p+1)/(p-1)}(p+1)^{q/(p-1)}L_{p,0}^{(2q-p+1)/(p-1)}
L_{p,q}
\right\}^{-1/(q-p+1)}
W_p(x),
\nonumber
\end{eqnarray}
\begin{eqnarray}
\lambda(\alpha) &=& (p+1)L_{p,0}L_{p,q}\alpha^{q-p+1}.
\end{eqnarray}
}

\noindent
{\bf Theorem 1.3.} {\it Let $f(x) = (1-x)^2$.  Assume that $q = m(p+1)$ or $q = m(p+1)+p$, where $m \in \mathbb{N}$. Then  

\noindent
(i)
\begin{eqnarray}
u_\lambda(x) &=& \left(\frac{\lambda}{M_{2,q}}\right)^{1/(q-p+1)}W_p(x),
\\
\lambda(\alpha) &=& M_{2,q}\xi_p^{-(q-p+1)}\alpha^{q-p+1},
\end{eqnarray}
where
\begin{eqnarray}
M_{2,q} &:=& C_{0,m} + C_{0,1,m}S_{1,0} + C_{0,2,m}S_{2, 0}, \qquad (q = m(p+1)),
\\
M_{2,q} &:=& C_{1,m} + C_{1,1,m}S_{1,p} + C_{1,2,m}S_{2, p}, \qquad (q = m(p+1)+p),
\end{eqnarray} 
$C_{0,m}, C_{0,1,m}, C_{1,2,m}, C_{1,m}, C_{1,1,m}, C_{1,2,m}$ are constants which 
depend on $\xi_p$, and are obtained inductively.  Here, 
 \begin{eqnarray}
 S_{1,0} &=&\frac18, 
\\
S_{2,0} &=& \frac{1}{24},
\\
S_{1,p} &=&\xi_p,
\\
S_{2,p} 
&=& \xi_p - \sqrt{2(p+1)}\xi_p^{(3-p)/2}L_{p,1}.
\end{eqnarray}
\noindent
For the special case $m = 1$, the following (ii) and (iii) hold.  

\noindent
(ii) Let $q = p + 1$. Then
\begin{eqnarray}
M_{2,p+1} &=& R_{2,p+1} =  S_{0,q} - 2S_{1,q} + 2S_{2,q}
\\
&=& \sqrt{\frac{p+1}{2}}\xi_p^{(p+3)/2}L_{p,p+1} - \frac{1}{3(p+1)}\xi_p^{p+1} 
- \frac{p+1}{p+3}\sqrt{2(p+1)}\xi_p^{(5-p)/2}L_{p,1}.
\nonumber
\end{eqnarray}

\noindent
(iii) Let $q = 2p+1 (= p+1 + p)$. Then 
\begin{eqnarray}
M_{2, 2p+1} &=& R_{2,2p+1} =  S_{0,q} - 2S_{1,q} + 2S_{2,q}
\\
&=& \sqrt{\frac{p+1}{2}}\xi_p^{3(p+1)/2}L_{p, 2p+1} - \frac23\sqrt{2(p+1)}\xi_p^{(p+5)/2}
\left(\frac{1}{p+2}L_{p,p+2} + 2L_{p,1}\right).
\nonumber
\end{eqnarray}
}

\vspace{0.2cm}

Now we consider the case $p = q$.  

\vspace{0.2cm}

\noindent
{\bf Theorem 1.4.} {\it Let $p = q > 1$. 

\noindent
(i) Assume that $n \ge 2$. Then $u_\lambda$ is obtained 
inductively by using $S_{k,1}$ ($1 \le k \le n-2$). 

\noindent
(ii) Especially, let $n = 3$. Then 
\begin{eqnarray}
M_{3,p} = \sqrt{\frac{2}{p+1}}\xi^{(p+1)/2} - 3\sqrt{2(p+1)}
\xi_p^{(3-p)/2}L_{p,1}
\end{eqnarray}
and (1.12), (1.13) hold by replacing $M_{2,q}$ and $q$ with $M_{3,p} $ and $p$, respectively. 
}

\vspace{0.2cm}

\noindent
{\bf Remarks.} (i) The case $n = 1$ in Theorem 1.4 is contained in Theorem 1.1. 

\noindent
{(ii)} The novelty of Theorem 1.4 is to give a scheme to obtain 
$u_\lambda$ inductively. Unfortunately, the concrete value of $S_{k,1}$ for $k \ge 1$ is not able to calculate. We only find that  
$M_{n,p}$ is expressed by $S_{k,1}$ ($1 \le k \le n-2$). Exceptionally, we have (cf. (4.3) below)
\begin{eqnarray}
S_{0,1} = \sqrt{\frac{p+1}{2}}\xi_p^{(3-p)/2}L_{p,1}.
\end{eqnarray}
The remainder of this paper is organized as follows. In section 2, we first prove Theorem 1.1. 
Next, we explain the 
existence of $u_\lambda$ for $\lambda > 0$ and fundamental properties of $W_p$. In sections 3, 4 and 5, the proofs 
of Theorems 1.2, 1.3  and 1.4 will be given. The main tools of the proofs are time map argument and complicated 
direct calculation.

\section{Proof of Theorem 1.1 and Preliminaries} 

In what follows, 
we write $\xi = \xi_p$ for simplicity. 
In this section, we consider the case $f(x) = (1-x)^n$, where $n \in \mathbb{N}$. 
For a given $\lambda > 0$, 
we look for the solution $u_\lambda$ of the form $u_\lambda = tW_p$ ($t > 0$).

\vspace{0.2cm}

\noindent
{\bf Proof of Theorem 1.1}. (i) Let $p \not= q+1$.  
We look for the solution of (1.1) of the form $u_\lambda = t_\lambda W_p$, where $t_\lambda > 0$ 
is a suitable constant determined in (2.2) below. By (1.1), we have 
\begin{eqnarray}
-t_\lambda^{q+1}\left(\int_0^1 (1-x)^nW_p(x)^q dx\right)W_p''(x) 
= \lambda t_\lambda^p W_p(x)^p.
\end{eqnarray}
Then there exists a unique $t_\lambda > 0$ satisfying 
\begin{eqnarray}
t_\lambda^{} =\left(\frac{\lambda}{M_{n,q}}\right)^{1/(q-p+1)},
\end{eqnarray}
and we find that $u_\lambda = t_\lambda W_p$ satisfies (1.1). 
We next show that, if there exists a solution $u_\lambda$ of (1.1), then $u_\lambda = t_\lambda W_p$, 
where $t_\lambda$ is a constant given in (2.2).  Indeed, we put $\displaystyle{r_\lambda:=\left(\frac{Q_{n,q}}{\lambda}\right)^{1/(p-1)}}$ and 
$w_\lambda:= r_\lambda^{-1} u_\lambda$. Then we see from (1.1) and (1.10) that $w_\lambda$ satisfies (1.3). Namely, 
$w_\lambda = W_p$. Then by definition of $r_\lambda$ and (1.10), we obtain 
\begin{eqnarray}
\lambda r_\lambda^{p-1} &=& Q_{n,q} = \int_0^1 (1-x)^n u_\lambda(x)^qdx 
\\
&=& \int_0^1 (1-x)^n r_\lambda^q W_p(x)^qdx = r_\lambda^q M_{n,q}
\nonumber
\end{eqnarray}
By this and (2.2), we see that $\displaystyle{r_\lambda = \left(\frac{\lambda}{M_{n,q}}\right)^{1/(q-p+1)} = t_\lambda}$ and 
$u_\lambda = t_\lambda W_p$. Thus the proof is complete. \qed

\noindent
{(ii)} Now assume that $p = q+1$. 

\noindent
{(a)} 
Assume that $\lambda = M_{n,q}$. For $t > 0$, we put $U_\lambda := tW_p$ and substitute it into 
(1.1). Then we obtain 
\begin{eqnarray}
-t^{q+1}\left(\int_0^1 (1-x)^nW_p(x)^q dx\right)W_p''(x) 
= \lambda t^p W_p(x)^p.
\end{eqnarray}
Then we see that (2.4) drives us to (1.3). Therefore,  $U_\lambda = tW_p$ is a solution of (1.1).
Thus the proof of (ii) (a) is complete. \qed

\noindent
{(b)} 
The proof of (ii) (b) is the same as that of Theorem 1.1 (i). So we omit the proof. \qed

\noindent
{(c)} By (a) and (b) above, a solution $u_\lambda$ of (1.1) exists if and only if the equality 
$\lambda= M_{n,q}$ holds. Therefore, if $\lambda \not= M_{n,q}$, then (1.1) has no solutions. 
Thus the proof of Theorem 1.1 is complete. \qed

\vspace{0.2cm}

To calculate $M_{n,q}$, we need some fundamental properties of $W_p$. 
Since (1.3) is autonomous, we know from [5] that 
\begin{eqnarray}
W_p(x) &=& W_p(1-x), \quad 0 \le x \le \frac12,
\\
\xi & =& \Vert W_p\Vert_\infty = \max_{0\le x \le 1}W_p(x) = W_p\left(\frac12\right),
\\
W_p'(x) &>& 0, \quad 0 \le x < \frac12.
\end{eqnarray}
By (1.3), for $0 \le x \le 1$, we have 
\begin{eqnarray}
\{W_p''(x) + W_p(x)^p\}W_p'(x) = 0.
\end{eqnarray}
By this and (2.6), we have 
\begin{eqnarray}
\frac12W_p'(x)^2 + \frac{1}{p+1}W_p(x)^{p+1} = \mbox{constant} = \frac{1}{p+1}
W_p\left(\frac12\right)^{p+1} 
= \frac{1}{p+1}\xi^{p+1}.
\end{eqnarray}
By this and (2.8), for $0 \le x \le 1/2$, we have, using $\theta = \xi_p s$,
\begin{eqnarray}
W_p'(x) = \sqrt{\frac{2}{p+1}(\xi^{p+1} - W_p(x)^{p+1})}.
\end{eqnarray}
By (2.5)--(2.7), (2.10) and putting $W_p(x) = \xi s$, we have 
\begin{eqnarray}
\Vert W_p\Vert_q^q &=& 2\int_0^{1/2} W_p(x)^qdx 
\\
&=& 2\int_0^{1/2} W_p(x)^q \frac{W_p'(x)}
{\sqrt{\frac{2}{p+1}(\xi^{p+1} - W_p(x)^{p+1})}}dx
\nonumber
\\
&=& 2\sqrt{\frac{p+1}{2}}\xi^{(2q-p+1)/2}\int_0^{1} \frac{s^q}{\sqrt{1-s^{p+1}}}ds
\nonumber
\\
&=& 2\sqrt{\frac{p+1}{2}}\xi^{(2q-p+1)/2}L_{p,q}.
\nonumber
\end{eqnarray}

\section{Proof of Theorem 1.2}
Now we put $n = 1$ and consider the case $f(x) = 1-x$.  
By (1.6), (2.5), (2.11) and putting $s = 1-t$, we have 
\begin{eqnarray}
R_{1,q} &=& \int_0^1 sW_p(s)^qds 
\\
&=& \int_0^{1/2} sW_p(s)^qds +\int_{1/2}^1 sW_p(s)^qds
\nonumber
\\
&=& \int_0^{1/2} sW_p(s)^qds +\int_0^{1/2} (1-t)W_p(t)^qdt 
\nonumber
\\
&=&  \int_0^{1/2} W_p(s)^qds 
\nonumber
\\
&=& \sqrt{\frac{p+1}{2}}\xi^{(2q-p+1)/2}L_{p,q}.
\nonumber
\end{eqnarray}
By this, (1.5) and (2.11), we obtain 
\begin{eqnarray}
M_{1,q} &=& \Vert W_p\Vert_q^q - R_{1,q} = \sqrt{\frac{p+1}{2}}\xi^{(2q-p+1)/2}L_{p,q}.
\end{eqnarray}
By this and (2.2), we have 
\begin{eqnarray}
u_\lambda(x) &=& \left(\frac{\lambda}{M_{1,q}}\right)^{1/(q-p+1)}W_p(x)
\\
&=& \lambda^{1/(q-p+1)}\left\{\sqrt{\frac{p+1}{2}}\xi^{(2q-p+1)/2}L_{p,q}\right\}^{-1/(q-p+1)}
W_p(x).
\nonumber
\end{eqnarray}
This along with (1.8) implies (1.11). Now 
we put $x = 1/2$ in (3.3). Then we have 
\begin{eqnarray}
\alpha &=& \lambda^{1/(q-p+1)}
\left\{\sqrt{\frac{p+1}{2}}L_{p,q}\right\}^{-1/(q-p+1)}\xi^{-(p-1)/(2(q-p+1))}.
\end{eqnarray}
This along with (1.8) implies (1.12). Thus the proof of Theorem 1.2 is complete. \qed

\section{Proof of Theorem 1.3} 

In this section, let $n = 2$, namely, $f(x) = (1-x)^2$. As in section 3, 
we look for the solution of (1.1) of the form $u_\lambda = t_\lambda W_p$, where $t_\lambda > 0$ 
is a constant defined by (2.4). 
By (1.5), (1.6), (2.11) and (3.1), we have 
\begin{eqnarray}
M_{2,q} &=& \int_0^{1} W_p(x)^qdx -2\int_{0}^1 xW_p(x)^qdx + \int_0^1 x^2W_p(x)^qdx
\\
&=& \int_0^{1} x^2 W_p(x)^qdx = R_{2,q}.
\nonumber
\end{eqnarray}

\vspace{0.2cm}

\noindent
{\bf Lemma 4.1.} {\it Assume that $q = m(p+1)$ or $q = m(p+1)+ p$, where $m \in \mathbb{N}$. 
Then $R_{2,q}$ is explicitly determined inductively. 
}

\noindent
{\it Proof. } By (1.7), (2.5), (4.1) and putting $x = 1-t$, we have 
\begin{eqnarray}
R_{2,q} &=& 
\int_0^{1/2} x^2W_p(x)^qdx + \int_{1/2}^1 x^2W_p(x)^qdx
\\
&=& 
\int_0^{1/2} x^2W_p(x)^q dx + \int_0^{1/2} (1 - 2t + t^2)W_p(t)^qdt
\nonumber
\\
&=& \int_0^{1/2} W_p(x)^qdx - 2\int_0^{1/2} xW_p(x)^qdx + 2\int_0^{1/2} x^2W_(x)^qdx
\nonumber
\\
&=& S_{0,q} - 2S_{1,q} + 2S_{2,q}.
\nonumber
\end{eqnarray} 
By (2.6), we have $W_p'(1/2) = 0$. By this, (1.1), (1.3), (2.10), (2.11) and integration by parts, we have 
\begin{eqnarray}
S_{0,q} &=& \int_0^{1/2} W_p(x)^q dx = \sqrt{\frac{p+1}{2}}\xi^{(2q-p+1)/2}L_{p,q},
\\
S_{1,q} &=& \int_0^{1/2} xW_p(x)^{q-p}W_p(x)^p dx 
\\
&=& -\int_0^{1/2} xW_p(x)^{q-p}W_p''(x)dx
\nonumber
\\
&=& -\left[xW_p(x)^{q-p}W_p'(x)\right]_0^{1/2} + 
\int_0^{1/2} (xW_p(x)^{q-p})'W_p'(x)dx
\nonumber
\\
&=& \int_0^{1/2}\{W_p(x)^{q-p} + x(q-p)W_p(x)^{q-p-1}W_p'(x)\}W_p'(x)dx
\nonumber
\\
&=& \frac{1}{q-p+1}\left[W_p(x)^{q-p+1}\right]_0^{1/2} + (q-p)\int_0^{1/2}xW_p(x)^{q-p-1}(W_p'(x))^2dx
\nonumber
\\
&=& \frac{1}{q-p+1}\xi^{q-p+1} + (q-p)\int_0^{1/2} xW_p(x)^{q-p-1}\frac{2}{p+1}
(\xi^{p+1} - W_p(x)^{p+1})dx
\nonumber
\\
&=& \frac{1}{q-p+1}\xi^{q-p+1} + \frac{2(q-p)}{p+1}\xi^{p+1}\int_0^{1/2}xW_p(x)^{q-p-1}dx 
- \frac{2(q-p)}{p+1}S_{1,q}.
\nonumber
\end{eqnarray}
By this, we obtain 
\begin{eqnarray}
S_{1,q} 
&=& \frac{p+1}{2q-p+1}\left\{\frac{1}{q-p+1}\xi^{q-p+1} + \frac{2(q-p)}{p+1}
\xi^{p+1}S_{1, q-p-1}\right\}.
\end{eqnarray}
Similar to the argument to derive (4.5), we obtain 
\begin{eqnarray}
S_{2,q} &=& \int_0^{1/2} x^2W_p(x)^{q-p}W_p(x)^p dx 
\\
&=& -\int_0^{1/2} x^2W_p(x)^{q-p}W_p''(x)dx
\nonumber
\\
&=& -\left[x^2W_p(x)^{q-p}W_p'(x)\right]_0^{1/2} + 
\int_0^{1/2} (x^2W_p(x)^{q-p})'W_p'(x)dx
\nonumber
\\
&=& \int_0^{1/2}\left\{2xW_p(x)^{q-p} + x^2(q-p)W_p(x)^{q-p-1}W_p'(x)\right\}W_p'(x)dx
\nonumber
\\
&=& 2\int_0^{1/2} xW_p(x)^{q-p} W_p'(x)dx + (q-p)\int_0^{1/2}x^2W_p(x)^{q-p-1}(W_p'(x))^2dx
\nonumber
\\
&=& 2\int_0^{1/2} xW_p(x)^{q-p} W_p'(x)dx 
\nonumber
\\
&& + (q-p)\int_0^{1/2}x^2W_p(x)^{q-p-1}
\frac{2}{p+1}(\xi^{p+1} - W_p(x)^{p+1})dx
\nonumber
\\
&=& 
2\int_0^{1/2} xW_p(x)^{q-p} W_p'(x)dx + \frac{2(q-p)}{p+1}\xi^{p+1}\int_0^{1/2}x^2W_p(x)^{q-p-1}dx 
- \frac{2(q-p)}{p+1}S_{2,q}.
\nonumber
\end{eqnarray}
By this, we have 
\begin{eqnarray}
\frac{2q-p+1}{p+1}S_{2,q} &=& 2\int_0^{1/2} xW_p(x)^{q-p} W_p'(x)dx 
\\
&&+ \frac{2(q-p)}{p+1}\xi^{p+1}\int_0^{1/2}x^2W_p(x)^{q-p-1}dx 
\nonumber
\\
&:=& 2U_1 +  \frac{2(q-p)}{p+1}\xi^{p+1}S_{2, q-p-1}.
\nonumber
\end{eqnarray}
By (4.3) and integration by parts, we have 
\begin{eqnarray}
U_1 &=& \int_0^{1/2} x \left(\frac{1}{q-p+1}W_p(x)^{q-p+1}\right)'dx 
\\
&=& \left[ x\frac{1}{q-p+1}W_p(x)^{q-p+1}\right]_0^{1/2}
- \int_0^{1/2}\frac{1}{q-p+1}
W_p(x)^{q-p+1}dx
\nonumber
\\
&=& \frac{1}{2(q-p+1)}\xi^{q-p+1}-\frac{1}{q-p+1}S_{0,q-p+1}
\nonumber
\\
&=& \frac{1}{q-p+1}\left(\frac12\xi^{q-p+1}-\sqrt{\frac{p+1}{2}}\xi^{(2q-3p+3)/2}L_{p,q-p+1}\right).
\nonumber
\end{eqnarray}
By this and (4.7), we have 
\begin{eqnarray}
S_{2,q} &=& \frac{p+1}{2q-p+1}\left\{\frac{1}{q-p+1}\left(\xi^{q-p+1}-\sqrt{2(p+1)}\xi^{(2q-3p+3)/2}L_{p,q-p+1}\right) 
\right.
\\
&&\left.\qquad \qquad \qquad + \frac{2(q-p)}{p+1}\xi^{p+1}S_{2,q-p-1}
\right\}.
\nonumber
\end{eqnarray}
We repeat the calculation (4.2)--(4.9). Then we obtain 
$R_{2,q}$ for $q = m(p+1)$ and $q = m(p+1) + p$ inductively. 
Indeed, assume that $q = m(p+1)$ (resp. $q = m(p+1)+p$). Then by repeating $m$ times the argument above, we have 
\begin{eqnarray}
R_{2,q} &:=& C_{0,m} + C_{0,1,m}S_{1,0} + C_{0,2,m}S_{2, 0}, \qquad (q = m(p+1)),
\\
R_{2,q} &:=& C_{1,m} + C_{1,1,m}S_{1,p} + C_{1,2,m}S_{2, p}, \qquad (q = m(p+1)+p),
\end{eqnarray} 
where $C_{0,m}, C_{0,1,m}, C_{1,2,m}, C_{1,m}, C_{1,1,m}, C_{1,2,m}$ are explicit constants 
containing $\xi$, which are obtained  inductively. 
According to the case where $q = m(p+1)$, or $q = m(p+1)+p$, (4.10) and (4.11) are determined by 
 \begin{eqnarray}
 S_{1,0} &=& \int_0^{1/2} xdx = \frac18, 
\\
S_{2,0} &=& \int_0^{1/2} x^2 = \frac{1}{24},
\\
S_{0,p} &=& \int_0^{1/2} W_p(x)^pdx = -\int_0^{1/2} W_p''(x)dx
 \\
 &=& -\left[W_p'(x)\right]_0^{1/2} = W_p'(0) = \sqrt{\frac{2}{p+1}}\xi^{(p+1)/2},
 \nonumber
 \\
S_{1,p} &=& \int_0^{1/2} xW_p(x)^pdx = \int_0^{1/2} x(-W_p''(x))dx
\\
&=& \left[x(-W_p'(x))\right]_0^{1/2} + \int_0^{1/2} W_p'(x)dx
\nonumber
\\
&=& \left[W_p(x)\right]_0^{1/2} 
\nonumber
\\
&=& \xi,
\nonumber
\\
S_{2,p} &=& \int_0^{1/2} x^2W_p(x)^pdx = \int_0^{1/2} x^2(-W_p''(x))dx 
\\
&=& \left[-x^2 W_p'(x)\right]_0^{1/2} + \int_0^{1/2}2xW_p'(x)dx 
\nonumber
\\
&=& \left[2xW_p(x)\right]_0^{1/2} - 2\int_0^{1/2} W_p(x)dx 
\nonumber
\\
&=& \xi - \sqrt{2(p+1)}\xi^{(3-p)/2}L_{p,1}.
\nonumber
\end{eqnarray}
Since $M_{2,q} = R_{2,q}$ by (4.1), we see from (4.12)--(4.16) that $M_{2,q}$ is 
explicitly determined 
inductively. 
Thus the proof is complete. \qed

\vspace{0.2cm}

\noindent
{\bf Proof of Theorem 1.3.} (i) The proof of Theorem 1.3 follows from Lemma 4.1 immediately. 
 
\noindent
(ii) Let $q = p + 1$. Then by (4.3), (4.5), (4.9), (4.12), (4.13) and direct calculation, we have 
\begin{eqnarray}
S_{0,p+1} &=& \sqrt{\frac{p+1}{2}}\xi^{(p+3)/2}L_{p,p+1},
\\
S_{1,p+1} &=& \frac{p+1}{p+3}\left\{\frac12\xi^2 + \frac{1}{4(p+1)}\xi^{p+1}\right\},
\\
S_{2,p+1} &=& \frac{p+1}{p+3}\left\{
\frac12\left(\xi^2 - \sqrt{2(p+1)}\xi^{(5-p)/2}L_{p,2} \right)
+ \frac{1}{12(p+1)}\xi^{p+1}
\right\}.
\end{eqnarray}
By (4.17)--(4.19), we have 
\begin{eqnarray}
M_{2,p+1} &=& R_{2,p+1} =  S_{0,q} - 2S_{1,q} + 2S_{2,q}
\\
&=& \sqrt{\frac{p+1}{2}}\xi^{(p+3)/2}L_{p,p+1} - \frac{1}{3(p+3)}\xi^{p+1} 
- \frac{p+1}{p+3}\sqrt{2(p+1)}\xi^{(5-p)/2}L_{p,2}.
\nonumber
\end{eqnarray}

\noindent
(iii) Let $q = 2p+1 (= p+1 + p)$. Then by the similar calculation as above, we have 
\begin{eqnarray}
S_{0,2p+1} &=& \sqrt{\frac{p+1}{2}}\xi^{3(p+1)/2}L_{p, 2p+1},
\\
S_{1,2p+1} &=& \frac{2p+5}{3(p+2)}\xi^{p+2},
\\
S_{2,2p+1} &=& \frac{2p+5}{3(p+2)}\xi^{p+2} - \frac13\sqrt{2(p+1)}\xi^{(p+5)/2}
\left(\frac{1}{p+2}L_{p,p+2} + 2L_{p,1}\right).
\end{eqnarray}
By, (4.22)--(4.24), we have 
\begin{eqnarray}
M_{2, 2p+1} &=& R_{2,2p+1} =  S_{0,q} - 2S_{1,q} + 2S_{2,q}
\\
&=& \sqrt{\frac{p+1}{2}}\xi^{3(p+1)/2}L_{p, 2p+1} - \frac23\sqrt{2(p+1)}\xi^{(p+5)/2}
\left(\frac{1}{p+2}L_{p,p+2} + 2L_{p,1}\right).
\nonumber
\end{eqnarray}
Thus the proof of Theorem 1.3 is complete. \qed

\section{Proof of Theorem 1.4}

In this section we consider the case $p = q > 1$ and $f(x) = (1-x)^n$ for $n \in \mathbb{N}$ 
with $n \ge 2$. We show that we are able to obtain $M_{n,p}$ inductively by using the constants 
$S_{m,1}$ ($m \ge 0$).

\noindent
{\bf Proof of Theorem 1.4 (i)}. 

\noindent
{\it Case 1.} Assume that $n = 2k + 1$, where $k \in \mathbb{N}$. We put $t = 1-s$. By (2.6), 
we have 
\begin{eqnarray}
M_{2k+1,p} &=& \int_0^1 (1 - s)^{2k+1} W_p(s)^pds 
\\
&=& \int_0^{1/2} (1 - s)^{2k+1} W_p(s)^pds +\int_{1/2}^1(1-s)^{2k+1}W_p(s)^pds 
\nonumber
\\
&=& \sum_{r = 0}^{2k} (-1)^r{}_{2k+1}C_r \int_0^{1/2}s^r W_p(s)^pds 
- \int_0^{1/2} s^{2k+1}W_p(s)^pds +  \int_{0}^{1/2} t^{2k+1} W_p(t)^pdt 
\nonumber
\\
&=& \sum_{r = 0}^{2k} (-1)^r{}_{2k+1}C_r \int_0^{1/2}s^r W_p(s)^pds.
\nonumber
\end{eqnarray}
Therefore, $M_{2k+1,p}$ is obtained by $S_{r, p}$ ($0 \le r \le 2k$).  

\noindent
{\it Case 2.} Assume that $n = 2k$, where $k \in \mathbb{N}$. Then 
\begin{eqnarray}
M_{2k,p} &=& \int_0^1 (1 - s)^{2k} W_p(s)^pds 
\\
&=& \int_0^{1/2} (1 - s)^{2k} W_p(s)^pds +\int_{1/2}^1(1-s)^{2k+1}W_p(s)^pds 
\nonumber
\\
&=& \sum_{r = 0}^{2k-1} (-1)^r{}_{2k}C_r \int_0^{1/2}s^r W_p(s)^pds 
+ \int_0^{1/2} s^{2k}W_p(s)^pds +  \int_{1/2}^1 t^{2k} W_p(t)^pdt 
\nonumber
\\
&=& \sum_{r = 0}^{2k-1} (-1)^r{}_{2k}C_r \int_0^{1/2}s^r W_p(s)^pds 
+ 2\int_0^{1/2} s^{2k}W_p(s)^pds.
\nonumber
\end{eqnarray}
By (5.1) and (5.2), we find that $M_{2k,p}$ is obtained by $S_{r, p}$ ($0 \le r \le 2k$).  
By (2.7), for $r \ge 2$, we have 
\begin{eqnarray}
S_{r,p} &=& \int_0^{1/2} x^r W_p(x)^pdx 
\\
&=& \int_0^{1/2} x^r(-W_p(x)'')dx = \left[x^r(-W_p'(x))\right]_0^{1/2} 
+ \int_0^{1/2} rx^{r-1}W_p'(x)dx
\nonumber
\\
&=& r\left[x^{r-1}W_p(x)\right]_0^{1/2} - r(r-1)\int_0^{1/2}x^{r-2}W_p(x)dx
\nonumber
\\
&=& r\left(\frac12\right)^{r-1}\xi - r(r-1)S_{r-2, 1}.
\nonumber
\end{eqnarray}
By this, we see that $S_{r,p}$ ($r \ge 2$) is represented by $S_{r-2, 1}$. 
By (5.2) and (5.3), we have 
\begin{eqnarray}
M_{2k+1,p} &=& \sum_{r = 0}^{2k} (-1)^r{}_{2k+1}C_r 
\left\{ r\left(\frac12\right)^{r-1}\xi - r(r-1)S_{r-2,1}\right\},
\\
M_{2k,p} &=& \sum_{r = 0}^{2k-1} (-1)^r{}_{2k}C_r  
\left\{r\left(\frac12\right)^{r-1}\xi - r(r-1)S_{r-2,1}\right\} 
\\
&&
+ 2\left\{2k\left(\frac12\right)^{2k-1}\xi - 2k(2k-1)S_{2k-2,1}\right\}.
\nonumber
\end{eqnarray}
By (5.4) and (5.5), we obtain our conclusion. Thus the proof is complete. \qed 

\noindent
{\bf Proof of Theorem 1.4 (ii)}. 
By (4.14), (4.15) and (4.16) and putting $t = 1-x$, we have 
\begin{eqnarray}
M_{3,p} &=& \int_0^1 (1-x)^3W_p(x)^pdx 
\\
&=&  \int_0^{1/2} (1-x)^3W_p(x)^pdx +  \int_{1/2}^1 (1-x)^3W_p(x)^pdx 
\nonumber
\\
&=&  \int_0^{1/2} (1-3x+3x^3-x^3)W_p(x)^pdx + \int_0^{1/2} t^3W_p(t)dt
\nonumber
\\
&=& \int_0^{1/2} (1-3x+3x^3)W_p(x)^pdx = S_{0,p} - 3S_{1,p} + 3S_{2,p} 
\nonumber
\\
&=& \sqrt{\frac{2}{p+1}}\xi^{(p+1)/2} - 3\sqrt{2(p+1)}\xi^{(3-p)/2}L_{p,1}.
\nonumber
\end{eqnarray}
Thus the proof is complete. \qed

\end{document}